\documentclass[12pt,amsfonts]{article}
\usepackage[margin=1in]{geometry}  
\usepackage{graphicx,color}              
\usepackage{amsmath}               
\usepackage{amsfonts}              
\usepackage{amsthm} 

\usepackage{tikz}

\def\nd{\noindent}
\def\thend{\rule{3mm}{3mm}}
\def\R{\mathbb{R}}
\def\N{\mathbb{N}}
\newtheorem{claim}{Claim}[section]
\newtheorem{thm}{Theorem}[section]

\newtheorem{lem}{Lemma}[section]

\newtheorem{deff}{Definition}[section]

\numberwithin{equation}{section}

\begin{document}

\title{Existence, regularity and  concentration phenomenon  of nontrivial solitary waves for a class of Generalized  Kadomtsev-Petviashvili (GKP) equation in $\mathbb{R}^2$  }
\author{\sf Claudianor O. Alves\thanks{Research of C. O. Alves partially supported by  CNPq/Brazil 304036/2013-7  and INCTMAT/CNPq/Brazil.} 
\\
\small{Universidade Federal de Campina Grande, }\\
\small{Unidade Acad\^emica de Matem\'atica } \\
\small{CEP: 58429-900 - Campina Grande-PB, Brazil}\\
\small{ e-mail: coalves@dme.ufcg.edu}\\
\vspace{1mm}\\
{\sf Ol\'impio H. Miyagaki\thanks{ O.H.M. was partially supported by INCTMAT/CNPq/Brazil, CNPq/Brazil 304015/2014-8.
}} \\
{ \small Departamento de Matem\'atica}\\
{\small Universidade Federal de  Juiz de Fora}\\
{ \small CEP:  36036-330 - Juiz de Fora - MG, Brazil}\\
 {\small e-mail: ohmiyagaki@gmail.com}\vspace{1mm} \\
}
\date{}
\maketitle
\begin{abstract}
In this paper we establish some results concerning the existence, regularity and concentration phenomenon of nontrivial solitary waves for a  Generalized Kadomtsev-Petviashvili (GKP) equation in $\mathbb{R}^2.$  Variational methods are used to get an existence result and to study the concentration phenomenon, while the regularity is  more delicate because we are leading with functions  in  an anisotropic Sobolev space. 
\end{abstract}

\vspace{0.5 cm}


\section{Introduction}

Consider the  Generalized Kadomtsev Petviashvili (GKP)  equation  of the form

\begin{equation}\label{eq1}
\left\{ \begin{array}{l}  u_t + V_x(x,y) {h(u)}+ V(x,y) h'(u)  u_x + u_{xxx} +\beta v_y=0  \\
v_x= u_y, \end{array} \right. 
\end{equation}
where $u=u(t,x,y)$ with $(t,x,y) \in \R^{+}\times \R \times \R,$  $h:\R \rightarrow \R$ and $V:\R^2 \rightarrow \R$ are smooth functions. A solitary wave of the (\ref{eq1}) is a solution of the form $u(t,x,y)=u(x-\tau t, y),$ with $\tau>0$. Hence, the function $u$ must satisfy the problem
\begin{equation}\label{eq2}
\left\{ \begin{array}{l}  -\tau u_x +  V_x(x,y) {h(u)}+ V(x,y) h'(u)  u_x + u_{xxx} +\beta v_y=0  \\
v_x= u_y.  \end{array} \right. 
\end{equation}
In the sequel, we will treat the case $\beta =-1$ and  $\tau=1.$  

By a simple calculus, it is easy to see that the above equation becomes
\begin{equation}\label{eq3}
  -u_x +  ( V(x,y) h(u))_x + u_{xxx} - D^{-1}_{x} u_{yy}=0  
\end{equation}
or equivalently
\begin{equation}\label{eq4}
  (-u_{xx} -  V(x,y) {h(u)}+ u + D^{-2}_{x} u_{yy})_x=0,
\end{equation}
where $D^{-1}$ denotes the following operator
$$
D^{-1}_{x}g(x,y)=\int_{-\infty}^{x} g(s,y)ds.      
$$

The equation (\ref{eq1}) is a  two-dimensional  Korteweg-de Vries type equation, which  is a model for long dispersive waves, essentially unidimensional, but having small transverse effects, see  \cite{GKP}. For the Cauchy problem associated with  equation  (\ref{eq1}) we would like to cite, e.g. \cite{Bourgain, Fam, IsazaM, WangAS} and the survey \cite{Saut}.  Recently, another interesting question is studied, namely, the existence and multiplicity of solitary waves to equation (\ref{eq1}). The pioneering work is due to De Bouard and Saut in \cite{SautB1,SautB2}, they treated a nonlinearity $h(s)=|s|^{p}s$ with $p=\frac{m}{n}$,$1 \leq  p < 4 ,$ if $N=2,$ and $1 \leq p< 4/3$ if $N=3,$ also,  $m$ and $n$ relatively prime, and $n$ is odd. In the mentioned paper, De Bouard and Saut obtained existence results by combining minimization with concentration compactness  theorem \cite{Lions}. For the regularity of the solutions they assumed $ p=2,3,4$ if $N=2,$ and $p=2$ if $N=3.$ In \cite{Willem} and \cite{WangWillem} a class of GKP problems were considered with an autonomous continuous nonlinearity $h$ in $N=2,$  and they have been proved the existence and multiplicity results, respectively. Their results were obtained by applying the mountain pass theorem \cite{AR} and Lusternik-Schnirelman theory, respectively.  In \cite{Liang}, Liang has proved the existence of solution for a class considered the GKP problem which involves  a non autonomous  continuous function with $N \geq 2,$ while \cite{Xuan} treated the autonomous case in higher dimension. We  recall that in the four above papers, the regularity of the solutions have not been treated.

Since a remarkable work by Rabinowitz in \cite{Rabinowitz}, also by \cite{W},  the existence and concentration of solutions for nonlinear Schr\"odinger equations, of the form

$$
\ \  \left\{
\begin{array}{l}
 - \epsilon^{2}\Delta u + V(z)u=f(u)
\ \ \mbox{in} \ \ \R^N (N \geq 2)
\\
u \in H^1(\R^N),
\end{array}
\right.
\eqno{(P_{\epsilon})}
$$
have been extensively
studied not only improving hypotheses on $f$, for instance by   \cite{AOS,Alves10,BPW,BW0,Pino,OS,FW,Oh1}  and  references therein,  but also treating other kind of the operator, we would like to cite \cite{AM} for the fractional laplacian, \cite{AG2006} for the p-laplacian and their  references.

Motivated by the above results, in the present paper we are going to study  the existence, regularity and concentration phenomenon of solitary waves for (\ref{eq1}), more precisely we consider the following equation  
\begin{equation}\label{eq6}
  -u_x + \left(V(\epsilon x,\epsilon y) h(u)\right)_x + u_{xxx} - D^{-1}_{x} u_{yy}=0, 
\end{equation}
or equivalently
\begin{equation}\label{eq7}
  \left(-u_{xx} - V(\epsilon x, \epsilon y) h(u) + u + D^{-2}_{x} u_{yy}\right)_x=0,
\end{equation}
with $ (x,y)\in \R^2$ and $\epsilon >0$.

We are going to assume, a similar set of the hypotheses on $h$ as that used in  \cite{Willem, WangWillem,Xuan}, and in $V$ as it was  imposed in \cite{Liang}. For $h$ we assume:

\begin{description}
\item [$(h_1)$]   $ \quad h \in C^2(\R^2)$  with $h(0)=h'(0)=0.$\\
\item [$(h_2)$]  \quad  There exists $C>0$  and $p \in (1,4)$ such that
$$
|h''(t)|\leq C |t|^{p-1}, \ \forall t \in \R.
$$
\item [$(h_3)$] There exists $\theta > 2$ such that
$$
0 < \theta H(t)\leq h(t)t, \quad \forall t \in \R\setminus\{0\} \quad \quad \mbox{where} \quad H(t)=\int_{0}^{t}h(r)dr.
$$

\item [$(h_4)$] \quad $\frac{h(t)}{|t|}$ \   is strictly increasing in $\R\setminus\{0\}.$
\end{description}

For $V$ we assume:
\begin{description}
\item [($V_1)$] $\displaystyle V \in L^{\infty}(\R^2),  \ V \geq 0,\  D^{\alpha}V \in L^{\infty}(\R^2), \mbox{for all } \ \alpha \in \mathbb{Z}^2,$  with $0\leq |\alpha|\leq 2,$
\item [($V_2)$ ] $\displaystyle \limsup_{|(x,y)|\rightarrow \infty} V(x,y)< \sup_{(x,y)\in \R^2} V(x,y).$
\end{description}

We will establish the following result

\begin{thm}\label{existencia}  Suppose $(h_1)-(h_4)$ and $(V_1)-(V_2)$ hold. The problem (\ref{eq6}) possesses at least one positive solution for $\epsilon$ small enough. Moreover, if $u_\epsilon$ denotes one of these solutions and $q_\epsilon$ is a global maximum point of $|u_\epsilon|$, we have that
$$
\lim_{\epsilon \to 0}V(\epsilon q_\epsilon)=\max_{(x,y) \in \mathbb{R}^2}V(x,y).
$$
\end{thm}

Here, we would like point out that for existence of solution it is enough to consider $p \in (0,4)$ in $(h_2)$. The restriction $p \in (1,4)$ is due to a technical difficulty to regularize the solutions of (\ref{eq6}). 

In the present paper, the motivation for using the term "{\it concentration phenomenon}" for family $u_\epsilon$ comes from the following fact: If we consider the family
$$
\zeta_\epsilon(x,t)=u_\epsilon(x/\epsilon, t/\epsilon)
$$
we will get a solution for the following class of problems
\begin{equation}\label{eq2}
\left\{ \begin{array}{l}  
-\epsilon u_x + \epsilon V_x {h(u)}+ \epsilon V h'(u)  u_x + \epsilon^{3} u_{xxx} - \epsilon v_y=0  \\
v_x= u_y. \end{array} \right. 
\end{equation}
In this case, for $q_\epsilon$ defined above,  if $\xi_\epsilon$ denotes a global maximum point of $|\zeta_\epsilon|$, we must have
$$
\xi_\epsilon= \epsilon q_\epsilon.
$$
As a byproduct of the proof of Theorem \ref{existencia}, we are able to prove that  for any sequence $\epsilon_n \to 0$, we have that
$$
\lim_{n \to +\infty}\xi_{\epsilon_n}=y \in \mathcal{V},
$$
where $\mathcal{V}=\{z\in \R^2: \ V(z)=\displaystyle \max_{(x,y) \in \mathbb{R}^2}V(x,y)\}.$ Then, the maximum point of $\zeta_\epsilon$ are concentrated near of $\mathcal{V}$  for $\epsilon$ small enough.

The plan of the paper is as follows: In Section 2, we will fix some notations and prove the existence of solution for $\epsilon$ small enough. In Section 3 we study the regularity of solutions, while in  Section 4 we study the concentration phenomenon. 

\vspace{0.5 cm}

\noindent {\bf Notations:} \, Throughout the paper, unless explicitly stated, the symbol $C$ will always denote a generic positive constant,
which may vary from line to line. The symbols \lq\lq $\rightarrow$ \rq\rq and \lq\lq $\rightharpoonup$ \rq\rq denote, respectively, strong and weak convergence, and  all the convergences involving sequences in $n\in \N$ are as $n\rightarrow \infty.$

\section{Notations and definitions}

Since we intend to use variational methods to prove our main result, we need to fix some notations and definitions. To begin with, we introduce the following function space 
\begin{deff} On $Y=\{ g_x: g \in C_{0}^{\infty}(\R^2)\}$ define the inner product  
\begin{equation}\label{produtointerno}
(u,v)=\int_{\R^2} ( u_x v_x +D^{-1}_{x}u_y D^{-1}_{x} v_y + uv)  dx dy
\end{equation}
with corresponding norm 
\begin{equation}\label{norm}
\|u\|=\left(\int_{\R^2} ( u^{2}_x +(D^{-1}_{x}u_y)^2  + u^2)  dx dy \right)^{\frac{1}{2}}.
\end{equation}
We say that $u:\R^2 \rightarrow \R$ belongs to $X$  if there exists  a sequence $(u_n)\subset  Y$ such that
\begin{description}
\item[a)] $ u_n \rightarrow u$ a. e. on $\R^2$,
\item[b)] $\|u_j - u_k\|\rightarrow 0, \  j, k \rightarrow \infty.$
\end{description}
\end{deff}

The space $X$  endowed with inner product and norm given above is a Hilbert space. Moreover, we have the following continuous embeddings whose proof can be found in \cite[Theorem 15.7  p. 323]{Besov} and \cite[Lemma 2.1]{Liang}
\begin{equation} \label{continuous}  
X \hookrightarrow L^{q}(\R^2), \quad 1\leq q \leq 6.
\end{equation}
Related to compact embeddings, De  Bouard and Saut in \cite[Remark 1.1]{SautB1} have proved that the embeddings below 

\begin{equation}\label{compact}  
X \hookrightarrow L^{q}_{loc}(\R^2), \quad \mbox{for} \quad 1\leq q < 6,
\end{equation}
are compacts. For higher dimension see \cite[Lemma 2.4]{Xuan}.

\subsection{The energy functional}
Associated with equation (\ref{eq6}) we have the energy functional $I_{\epsilon}:X \longrightarrow \R$ given by
$$
I_{\epsilon}(u)=\frac{1}{2}\|u\|^2-\int_{\R^2}V(\epsilon x, \epsilon y) H(u) dx dy.
$$
We are going to assume that 
$$
\limsup_{|(x,y)|\rightarrow \infty} V(x,y) >0; \eqno{(V_0)}
$$
since the case
$$
\limsup_{|(x,y)|\rightarrow \infty} V(x,y) =0
$$
leads to 
\begin{equation}\label{1}
\lim_{|(x,y)|\rightarrow \infty} V(x,y) =0.
\end{equation}

When (\ref{1}) occurs, the functional $\Psi:X \longrightarrow \R$ given by
$$
\Psi(u)=\int_{\R^2}V(\epsilon x, \epsilon y) H(u) dx dy
$$
is weakly continuous, that is,
\begin{equation}\label{2}  
{If}\ u_n \rightharpoonup u \ \mbox{weakly in} \ X,\ \mbox{then}\  \Psi(u_n)\rightarrow \Psi(u) \  \mbox{in}\ \R.
\end{equation}
Moreover, we also have 
\begin{equation}\label{3}  {If}\ u_n \rightharpoonup u \ \mbox{weakly in} \ X,\ \mbox{then}\  \Psi'(u_n)\rightarrow \Psi'(u) \  \mbox{in}\ X'.
\end{equation}

From $(h_1)-(h_3)$, (\ref{1}),(\ref{2}) and (\ref{3}), it follows that $I_\epsilon$  verifies the Mountain Pass geometric conditions and Palais-Smale (PS) condition. Then,  applying the Ambrosetti-Rabinowitz mountain pass theorem \cite{AR} there exists a critical point $u_{\epsilon} \in X$ with
$$
0< I_{\epsilon}(u_{\epsilon})=c_{\epsilon},
$$
where $c_{\epsilon}$ is the mountain pass level associated with $I_\epsilon$. 

 \subsection{The existence of solution with $(V_0)$}
By standard arguments, $I_{\epsilon}$ verifies the Mountain Pass geometric conditions for all $\epsilon>0,$ then there exists  a $(PS)_{c_{\epsilon}}$ sequence $(u_n)\subset X$, that is, 
$$
I_{\epsilon} (u_n)\to c_{\epsilon}\quad \mbox{and} \quad \ I'_{\epsilon}(u_n) \to 0 \quad \mbox{as} \ \epsilon \to 0,
$$
where  $c_{\epsilon}$ is the mountain pass level associated with $I_\epsilon$. Moreover, $(u_n)$ is  bounded in $X$  and there exists $u_\epsilon \in X$ such that $u_n \rightharpoonup u_{\epsilon}$ in $X.$ As $I'_\epsilon(u_n)v=o_n(1)$ for each $v \in X$, we derive that 
$$
I'_\epsilon(u_\epsilon)v=0 \quad \forall v \in X,
$$
and so, 
$$ 
I'_{\epsilon}(u_\epsilon)=0.
$$
In the sequel, we will show that $u_\epsilon \neq 0$ for $\epsilon$ small enough.

In what follows, $c_0$ and $c_\infty$  denote the mountain pass levels associated with the functionals
$$
I_{0}(u)=\frac{1}{2}\|u\|^2-\int_{\R^2}V_0 H(u) dx dy
$$
and
$$
I_{\infty}(u)=\frac{1}{2}\|u\|^2-\int_{\R^2}V_\infty  H(u) dx dy,
$$
respectively, where
\begin{equation} \label{Vinfty}
V_0=V(0,0)\quad \mbox{and}\quad V_\infty=\limsup_{|(x,y)|\rightarrow +\infty}V(x,y).
\end{equation}
Without lost of generality, we assume that 
$$
V_0=\max_{(x,y)\in \R^2}V(x,y).
$$
Then by $(V_2)$, 
$$
V_0 > V_\infty,
$$
from where it follows that
\begin{equation}\label{Delta}
c_0< c_\infty.
\end{equation}

\begin{lem}\label{lema1} Suppose that
$$\limsup_{\epsilon \to 0} c_{\epsilon}< c_\infty. \leqno{(H)}$$
Then, $u_\epsilon \neq 0$ for all $\epsilon$  sufficiently small.
\end{lem}
\noindent{\bf Proof.}  Suppose by contradiction $u_\epsilon=0$ for  a $\epsilon >0$ fixed. By using Lions' Lemma version for $X$ found in \cite{Willem}, there exist $(y_n)\subset \R^2$ and $ R, \eta >0$  such that
\begin{equation}\label{lions}
\lim_{n\in \mathbb{N}}\int_{B_{R}(y_n)}|u_n|^2 dx dy\geq \eta >0.
\end{equation}
Here $B_R (a)$ denotes an open ball centered at  $a$ with radius $R.$

Now define the translated sequence $w_n(x,y)=u_n((x,y)+y_n)$ and note that it is bounded in $X$. Thus,  up to a subsequence,
$$
w_n \rightharpoonup w \ \mbox{in } X \ \mbox{and}\ w_n \rightarrow w \ \mbox{in} \ L^{2}_{loc}(\R^2).
$$
This together with (\ref{lions}) implies that $w\neq 0.$ Considering the test function $v_n(x,y)=w((x,y)-y_n),$ we infer that $\|v_n\|=\|w\|$ for all $n \in \mathbb{N}$. Hence, $(v_n)$ is also bounded in $X$, and so,
$$
I'_{\epsilon}(u_n)v_n=o_n(1).
$$
After changing variable, we have
\begin{equation}\label{4}
(w_n,w)=\int_{\R^2}V(\epsilon(x,y)+ y_n)h(w_n) w  dx dy.
\end{equation}
Notice that the sequence 
$$
f_n(x,y)=V(\epsilon(x,y)+ y_n)h(w_n) w
$$
is bounded from above by sequence 
$$ 
g_n(x,y)=V_0 |h(w_n)| |w|,
$$
which converges in $L^{1}(\R^2)$ to $ g(x,y)=V_0 |h(w)||w|.$ From this, the reverse Fatou's Lemma gives 
\begin{equation}\label{5} \limsup_{n \to +\infty } \int_{\R^2} V(\epsilon(x,y)+ y_n)h(w_n) w dx dy\leq  \int_{\R^2} \limsup_{n \to +\infty}  V(\epsilon(x,y)+ y_n)h(w) w dx dy.
\end{equation}

The above inequality helps us to prove the following claim 

\begin{claim} \label{CL1} \, The sequence $(y_n)$ is bounded in  $\mathbb{R}^2$ for $\epsilon$ small enough.

\end{claim}
Indeed, suppose that $(y_n)$  possesses a subsequence, still denoted by $(y_n)$, such that
$$
|y_n|\to +\infty.
$$
Thereby, by (\ref{4}),  
\begin{equation}\label{6} \limsup_{n\to +\infty} \int_{\R^2} V(\epsilon(x,y)+ y_n) h(w_n) w dx dy\leq  \int_{\R^2}  V_{\infty} h(w) w dx dy.
\end{equation}
Therefore, from (\ref{4}), (\ref{5}) and (\ref{6}), 
$$
\|w\|^2\leq  \int_{\R^2}  V_{\infty}h(w) w dx dy,
$$
which yields 
\begin{equation}\label{7} 
I'_{\infty}(w)w\leq 0.
\end{equation}
Let $t \in (0, +\infty)$ be a number such that
$$
I_{\infty}(tw)=\max_{s \geq 0}I_{\infty}(sw).
$$
From (\ref{7}), we infer that $t \in (0,1].$ Recalling that $c_\infty$ can be characterized by infimun on Nehari manifold associated with $I_\infty$ (see \cite{Willem}), it follows that
$$
c_\infty \leq  I_\infty(tw)=I_\infty(tw)- \frac{1}{2}I'_{\infty}(tw)tw=\int_{\mathbb{R}^2}\frac{V_\infty}{2}(h(tw)(tw)-2H(tw).
$$
By $(h_4)$, the function $f(s)=h(s)s-2H(s)$ is increasing for $s >0 $ and decreasing for $s < 0$, we find  
\begin{eqnarray}\label{8}
c_\infty& \leq &\int_{\mathbb{R}^2}\frac{V_\infty}{2}(h(w)(w)-2H(w))dxdy\\
&\leq & \liminf_{n \to +\infty}\int_{\mathbb{R}^2}\frac{V(\epsilon x,\epsilon  y)}{2}(h(u_n)(u_n)-2H(u_n))dxdy. \nonumber\\
&= & \lim_{n \to +\infty}(I_\epsilon(u_n)-\frac{1}{2}I'_\epsilon(u_n)(u_n)) \nonumber\\
&= & \lim_{n \to +\infty}I_\epsilon(u_n)=c_\epsilon, \nonumber\\
\end{eqnarray}
that is, $c_\epsilon \geq c_\infty$. On the other hand, condition $(H)$ implies that there is $\epsilon_0>0$ such that
$$
c_\epsilon < c_\infty, \quad \forall \epsilon \in (0, \epsilon_0).
$$
From the above analysis,  $(y_n)$ must be a bounded sequence for $\epsilon \in(0, \epsilon_0).$

From Claim \ref{CL1}, as $(y_n)$ is bounded, there exists  $\widehat{R}>0$  such that
$$
\int_{ B_{\widehat{R}}(0)}|u_n|^2 dx dy=\int_{ B_{R}(y_n)}|u_n|^2 dx dy\geq \eta>0.
$$
Then, applying the compact embedding (\ref{compact}),  
$$
\int_{ B_{\widehat{R}}(0)}|u_\epsilon|^2 dx dy\geq \eta>0
$$
from where it follows that $u_\epsilon \neq 0$ for $\epsilon \in (0, \epsilon_0)$.
\begin{flushright}$\Box$
\end{flushright}

The next lemma is very important in our approach, because it shows that condition $(H)$ holds for $\epsilon>0$ small enough.

\begin{lem}\label{lema2}
$$
\lim_{\epsilon \rightarrow 0} c_\epsilon=c_0.
$$
\end{lem}

\noindent {\bf Proof.} From the hypotheses  on $V$ 
$$
V(\epsilon x, \epsilon y)\leq V_0=V(0,0), \quad \forall x \in \mathbb{R}^N.
$$
Thus,
$$
c_\epsilon \geq c_0, \quad \forall \epsilon >0,
$$
leading to
\begin{equation}\label{10}
\liminf_{\epsilon \rightarrow 0} c_\epsilon\geq c_0.
\end{equation}

On the other hand, let $w_0$ be a  {\it ground state} solution associated with functional $I_0$, that is,
$$
I_0(w_0)=c_0 \quad \mbox{and}\quad I'_{0}(w_0)=0.
$$
Let $t_\epsilon >0$ be a number    such that
$$
I_{\epsilon}(t_\epsilon w_0)=\max_{s \geq 0}I_{\epsilon}(sw_0).
$$
This implies that
$$
\frac{d}{dt}I_{\epsilon}(t w_0)|_{t=t_\epsilon}=0
$$
or equivalently
$$
\|w_0\|^2=t^{-2}_{\epsilon}\int_{\R^2} V(\epsilon(x,y))h(t_\epsilon w_0)(t_\epsilon w_0) dx dy.
$$
Gathering $(h_4)$, the above identity and Fatou's Lemma, we conclude there is $t_{0}>0$ such that $t_\epsilon \to t_0>0.$ By Lebesgue's dominated convergence theorem, we find 
$$
\|w_0\|^2=t^{-2}_{0}\int_{\R^2} V_0 h(t_0 w_0)t_0 w_0 dx dy,
$$
that is,
$$
\frac{d}{dt}I_{0}(t w_0)|_{t=t_0}=0.
$$
But
$$
\frac{d}{dt}I_{0}(t w_0)|_{t=1}=0,
$$
then by uniqueness of the maximum point we deduce that $t_0=1,$ and so, 
\begin{equation}\label{11}
\lim_{\epsilon \rightarrow 0} t_\epsilon=1.
\end{equation}
Recalling the characterization of the mountain pass level $c_\epsilon,$ we have
$$
c_\epsilon\leq I_{\epsilon}(t_\epsilon w_0)=\frac{t^{2}_{\epsilon}}{2}\|w_0\|^2-\int_{\R^2}V(\epsilon(x,y))H(t_\epsilon w_0)  dx dy.
$$
Therefore,  
\begin{equation}\label{12}
\limsup_{\epsilon \rightarrow 0} c_\epsilon\leq I_0(w_0)=c_0.
\end{equation}
Gathering  (\ref{10}) and (\ref{12}), 
$$
\lim_{\epsilon \rightarrow 0} c_\epsilon=c_0.
$$
This completes the proof of Lemma \ref{lema2}.
\begin{flushright}$\Box$
\end{flushright}

Now, we are able to prove the existence of solution for $(P_\epsilon)$. 

\begin{thm}\label{teorema1}
There is $\epsilon^{*}>0$, such that the mountain pass level $c_\epsilon$ is a critical  value of $I_\epsilon$ for $\epsilon \in (0,\epsilon^*)$, that is, there exists $u_\epsilon \in X$ , so called {\it ground state solution}, such that
$$
I_\epsilon (u_\epsilon)=c_\epsilon \quad \mbox{and} \quad I'_\epsilon (u_\epsilon)=0.
$$
\end{thm}
\noindent {\bf Proof.} Since $c_0 < c_\infty$, by Lemma \ref{lema2} there is $\epsilon^*>0$ such that $c_\epsilon < c_\infty$ for all $\epsilon \in (0, \epsilon^*)$. Now, by applying Lemma \ref{lema1}  $u_\epsilon \neq 0$ for all $\epsilon \in (0, \epsilon^*)$. As $I'_\epsilon(u_\epsilon)=0$, it follows that
$$
u_\epsilon \in \mathcal{N}_\epsilon=\{ u \in X \setminus \{0\}:  \ I'_\epsilon (u)u=0\},
$$
and so,
\begin{equation}\label{13}c_\epsilon =\inf_{u \in \mathcal{N}_\epsilon}I_\epsilon (u)\leq I_\epsilon (u_\epsilon).
\end{equation}
On the other hand, the Fatou's Lemma leads to 
\begin{eqnarray}\label{14}
I_\epsilon (u_\epsilon)&=& I_\epsilon (u_\epsilon)-\frac{1}{\theta}I'_\epsilon (u_\epsilon)\nonumber\\
&\leq&\liminf_{n \to +\infty}( I_\epsilon (u_\epsilon)-\frac{1}{\theta}I'_\epsilon (u_\epsilon))\\
&=& \liminf_{n \to +\infty}  I_\epsilon (u_\epsilon)=c_\epsilon. \nonumber
\end{eqnarray}
From (\ref{13}) and (\ref{14}), 
$$
I_\epsilon (u_\epsilon)=c_\epsilon,
$$
that is, $c_\epsilon$ is a critical value for  $\epsilon \in (0, \epsilon^*).$ Moreover, $u_\epsilon$ is called of ground state solution for $(P_\epsilon)$.
\begin{flushright}$\Box$
\end{flushright}

\section{Regularity}

In this section we study the regularity of the solutions of (\ref{eq6}), because it is crucial to study the concentration phenomenon.  The regularity will get by using Fourier transform of a tempered distribution.  Next, we recall the definition  of the Fourier transform of a tempered distribution, for more details see \cite{folland}. For a distribution $f$ and multi-index  $\alpha$,  the derivative of $f$ is given by
$$
<\partial^\alpha f, \phi>=(-1)^{|\alpha|}<f, \partial^{\alpha}\phi>, \forall \phi \in C_{0}^{\infty}(\R^2).
$$

The Fourier transform $\widehat{f}$ of a tempered distribution is defined by
$$
<\widehat{f}, \phi>=<f,\widehat{\phi}>, \quad \forall \phi \in \mathcal{S} \quad  (\mbox{Schwartz space}),
$$
\noindent likewise, it is  defined the inverse Fourier transform of a tempered distribution  $f$, denoted by $f^{\vee}$, by
$$
f^{\vee}(x)=(2\pi)^{-2}\widehat{f}(-x).
$$
All the basic properties of the usual Fourier transform remain valid for Fourier transform of a tempered distribution, for instance, for all tempered distribution $f$ and multi-index $\alpha,$  we have
$$
\widehat{\partial^{\alpha}f(x)}=i^{|\alpha|}\xi^{\alpha}\widehat{f}(\xi) \ \ \mbox{and}\ \  \widehat{x^{\alpha} f(x)}=i^{|\alpha|}\partial^{\alpha}\widehat{f}(\xi).
$$

Now, we are able to state and prove the result below

\begin{thm}\label{thm4.1} Any solution $u$ of (\ref{eq6}) is continuous. Moreover,
	$$u \in W^{2,q'}(\R^2) \quad \mbox{with}\quad \left\{\begin{array}{rcl} q'=6/(p+1) &\mbox{if}& p\neq 3 \\
	  q'\in(1,3/2)& \mbox{if} & p= 3. \end{array} \right.$$ In addition
$$
u(x) \rightarrow 0 \quad \mbox{as}\quad |x|\rightarrow  \infty.
$$
\end{thm}

\noindent {\bf Proof.}\, We are going to adapt the proof of Theorem 4.1 given   in \cite{SautB1} for our case. Here the  standard procedure made to the laplacian operator does not work any longer, because  the symbol of the linear operator $ -\Delta + \partial^{4}_{x} $ is non isotropic. The proof is made bootstrapping, using a variant of the H\"ormander-Mikhlin multipliers theorem due to Lizorkin see
\cite[Corollary 1]{lizorkin}.

First of all, we must observe that if $u\in X$  is a solution of (\ref{eq7}), then it is a solution, in the distribution sense, of the problem  
\begin{equation}\label{eq9}
  -\Delta u + u_{xxxx} = g_{xx}, \quad  \mbox{in} \quad  \R^{2}
\end{equation}
where
$$
g(x,y)=-V(x,y) h(u(x,y)).
$$
Since the  functions involved in above equation (\ref{eq9}) belong to $L^s(\R^2)$ space for some $s>1,$ we can assume that they are tempered distribution in $\R^2.$

Applying the Fourier transfom in the equation (\ref{eq9}), in the sense of the tempered distribution,  we have
$$
<\widehat{-\Delta u}, \phi>+<\widehat{u_{xxxx}},\phi>=<\widehat{g _{xx}},\phi>, \forall \phi \in \mathcal{S}.
$$

By  the above mentioned properties, we obtain for each $\phi \in \mathcal{S}, $
$$
\begin{array}{c}
	<-\Delta u,\widehat{ \phi}>+<u_{xxxx},\widehat{\phi}>=<g _{xx},\widehat{\phi}>\\ \\
		< u,-\Delta(\widehat{ \phi})>+  <u,\widehat{\phi}_{xxxx}>= <g ,\widehat{\phi}_{xx}>\\ \\
		 < u,\widehat{ |(x,y)|^2 \phi}>+ <u,\widehat{x^4 \phi}>= -<g ,\widehat{x^2 \phi}>, \  (x,y) \in \R^2\\ \\
		 < \widehat{u}, |(\xi_1,\xi_2)|^2 \phi>+ <\widehat{u},\xi_{1}^{4}\phi>=-<\widehat{g} ,\xi_{1}^{2}\phi>, \ (\xi_1,\xi_2) \in \R^2 
\end{array}
$$
that is
$$
|\xi|^2 \widehat{u}(\xi)+|\xi_1|^4\widehat{u}(\xi)=-|\xi_1|^2\widehat{g}(\xi), \quad \xi=(\xi_1,\xi_2)\in \R^2,
$$
then
$$
\widehat{u}(\xi)=-(\frac{|\xi_1|^2}{|\xi|^2+|\xi_1|^4})\widehat{g}\equiv - \Phi_1(\xi)\widehat{g}.
$$
Hence
\begin{equation}\label{eq10}u=\overbrace{(-\Phi_1(\xi)\widehat{g})}^{\vee}.\end{equation}

By \cite[Corollary 1]{lizorkin}, $\Phi_1$ is a Fourier multiplier on $L^q(\R^2)$ for all $ q \in (1, \infty),$ that is, 
\begin{equation}\label{Multiplier}\|\overbrace{\Phi_1(\xi)\widehat{f}}^{\vee}\|_{L^q}\leq C \|f\|_{L^q}, \quad \forall f \in L^q(\R^2), \quad C \ \mbox{is independent on}\  f \mbox{and} \ \Phi_1 .
\end{equation}
Since $ u \in X$, from Theorem \ref{continuous}, we infer $g \in L^{6/(p+1)}(\R^2)$. Combining (\ref{eq10}) with   (\ref{Multiplier}), we conclude

\begin{equation}\label{regularu)}
u \in L^{6/(p+1)}(\R^2).
\end{equation}
Differentiating (\ref{eq9}) with respect to $x$, twice, we have

\begin{equation}\label{eq12}
  -\Delta ( u_{xx}) + (u_{xx})_{xxxx} =  g_{xxxx} , \quad  (x,y)\in \R \times \R.
\end{equation}
Applying, as above,  the Fourier transform in the equation (\ref{eq12}),  we get
$$|\xi|^2 \widehat{u_{xx}}(\xi)+|\xi_1|^4\widehat{u}_{xx}(\xi)=|\xi_1|^4\widehat{g}(\xi), \quad \xi=(\xi_1,\xi_2)\in \R^2,
$$
that is, 
$$
\widehat{u_{xx}}(\xi)=(\frac{|\xi_1|^4}{|\xi|^2+|\xi_1|^4})\widehat{g}\equiv \Phi_2(\xi)\widehat{g}
$$
or equivalently 
\begin{equation}\label{eq13} 
u_{xx}=\overbrace{(\Phi_2(\xi)\widehat{g})}^{\vee}.
\end{equation}

By \cite[Corollary 1]{lizorkin}, $\Phi_2$ is a Fourier multiplier on $L^q(\R^2)$ for all $ q =6/(p+1),$ that is, 

\begin{equation}\label{regularuxx}u_{xx} \in L^{6/(p+1)}(\R^2).
\end{equation}
Differentiating (\ref{eq9}) with respect to $y$,  we find 
\begin{equation}\label{eq16}
  -\Delta ( u_{y}) + (u_{y})_{xxxx} =  g_{xxy} , \quad  (x,y)\in \R \times \R.
\end{equation}
Similarly, applying the Fourier transform in the equation (\ref{eq16}), recalling that
$$
<g_{xxy},\widehat{\phi}>=-<g,\widehat{\phi}_{xxy}>=-i<g, \widehat{x^2 y \phi}>=-i<\widehat{g},\xi_{1}^{2} \xi_{2}\phi>,
$$
we have
$$
|\xi|^2 \widehat{u_{y}}(\xi)+|\xi_1|^4\widehat{u}_{y}(\xi)=-i |\xi_1|^2\xi_2\widehat{g}(\xi), \quad \xi=(\xi_1,\xi_2)\in \R^2,
$$
that is, 
$$
\widehat{u_{y}}(\xi)=(\frac{|\xi_1|^2 \xi_2}{|\xi|^2+|\xi_1|^4})\widehat{(-i)g}\equiv \Phi_3(\xi)\widehat{(-i)g}.
$$
Hence
\begin{equation}\label{eq17} u_{y}=\overbrace{(\Phi_3(\xi)\widehat{(-i)g})}^{\vee}.\end{equation}

By \cite[Corollary 1]{lizorkin}, $\Phi_3$ is a Fourier multiplier on $L^q(\R^2)$ for all $ q =6/(p+1),$ that is, 

\begin{equation}\label{regularuy}u_{y} \in L^{6/(p+1)}(\R^2).
\end{equation}
We claim
\begin{equation}\label{EQ1}
u_x  \in  L^{12/(p+1)}(\R^2).
\end{equation}

\noindent
{\bf Verification.}  Notice $ u \in W^{\overrightarrow{l}}_{\overrightarrow{p}}(\R^2),$ with $\overrightarrow{l}=(2,0)$ and $ \overrightarrow{p}=(6/(p+1), 6/(p+1)). $  Since $u, u_{xx}\in  L^{6/(p+1)}(\R^2),$ we are going to apply  \cite[Thm 10.2]{Besov} with $\alpha=(1,0),$  $\overrightarrow{q}=(12/(p+1), 12/(p+1)),$  $ l_1=(2,0)$ and $ l_2=(0,0),$
in this case the line containing the points $ (2,0)$ and $(0,0)$ is just $y=0.$ 
The point $w=\alpha + \frac{1}{\overrightarrow{p}}-\frac{1}{\overrightarrow{q}}=(1,0)+((p+1)/6, (p+1)/6)-((p+1)/12, (p+1)/12)=((p+13)/12, (p+1)/12)$
does not lie on the above line. Then, there exist  positive constants $C_1, C_2>0$ such that 
$$
|u_x|_{L^{12/(p+1)}}\leq C_1(|u_{xx}|_{L^{6/(p+1)}} + |u|_{L^{6/(p+1)}}) + C_2 |u|_{L^{6/(p+1)}}.
$$
We recall that, up to now,
\begin{equation}\label{Passo1a}
 u, u_{xx}, u_y \in  L^{6/(p+1)}(\R^2) \ \mbox{and} \ u_x  \in  L^{12/(p+1)}(\R^2), \ \forall p \in [1,4).
\end{equation}
We claim
\begin{equation}\label{Passo1b}
u \in L^{r}(\R^2)\quad \mbox{for}\quad  \left\{ \begin{array}{rcl}   r \in [6/(p+1), \infty] & \mbox{if}& p \neq 3 \\  r \in (1, \infty) &\mbox{if}&  p=3. \end{array}\right.
\end{equation}

\noindent {\bf Verification for $p\neq 3$:} \,  Notice $ u \in W^{\overrightarrow{l}}_{\overrightarrow{p}}(\R^2),$ with $\overrightarrow{l}=(2,1)$ and $ \overrightarrow{p}=(6/(p+1), 6/(p+1)). $  Since $u_y, u_{xx}\in  L^{6/(p+1)}(\R^2),$ we are going to apply  \cite[Thm 10.2]{Besov} with $\alpha=(0,0),$  $\overrightarrow{q}=(\infty,\infty),$  $ l_1=(2,0)$ and $ l_2=(0,1).$ In this case the line containing the points $ (2,0)$ and $(0,1)$ is just $y=-\frac{x}{2}+1.$ The point $w=\alpha + \frac{1}{\overrightarrow{p}}-\frac{1}{\overrightarrow{q}}=(0,0)+((p+1)/6, (p+1)/6)-(0, 0)=((p+1)/6, (p+1)/6)$ does not  lie on the above line. By \cite[Thm 10.2]{Besov}, $ u \in L^{\infty}(\R^2).$ Thus, for $p \neq 3,$ we have $u\in L^{6/(p+1)}(\R^2)\cap L^{\infty}(\R^2),$ from where it follows that  
$$
u \in L^r(\R^2), \forall r\in [6/(p+1),\infty] \ \mbox{if} \ p \neq 3.
$$

\noindent {\bf Verification for $p= 3$:} \, Notice $ u \in W^{\overrightarrow{l}}_{\overrightarrow{p}}(\R^2),$ with $\overrightarrow{l}=(2,1)$ and $ \overrightarrow{p}=(3/2,3/2). $  Since $u_y, u_{xx}\in  L^{3/2}(\R^2),$ we are going to apply  \cite[Thm 10.2]{Besov} with $\alpha=(0,0),$  $\overrightarrow{q}=(r, r),$  $ l_1=(2,0)$ and $ l_2=(0,1).$ 
In this case the line containing the points $ (2,0)$ and $(0,1)$ is just $y=-\frac{x}{2}+1.$
The point $w=\alpha + \frac{1}{\overrightarrow{p}}-\frac{1}{\overrightarrow{q}}=(0,0)+(3/2,3/2)-(r, r)=(3/2 -r, 3/2- r)$ does not  lie on the above line, if $ r \in (1,\infty)$. Thereby, by  \cite[Thm 10.2]{Besov}, $ u \in L^{r}(\R^2)$ for all $r > 1.$

\vspace{0.5cm}

In the sequel, setting
$$
f=(V h(u))_{xx}, 
$$
we see that
$$
f= V_{xx} h(u)+2 V_x h'(u)u_x + V h''(u)u^{2}_{x}+ V h(u) u_{xx}.
$$

We claim
\begin{equation}\label{Pass1F}
f \in L^{q}(\R^2)\quad \mbox{for}\quad  \left\{ \begin{array}{rcl}   q =6/(p+1) & \mbox{if}& p \neq 3 \\  q \in (1, 3/2) &\mbox{if}&  p=3. \end{array}\right.
\end{equation}

\noindent {\bf Verification:} \, Since $V,V_x$ and $V_{xx}$ are assumed bounded, we can drop it. Let us analyze each of the terms of the $f.$ For first term, $(h_2)$ gives 
$$
|h(u)|\leq C |u|^{p+1}
$$
for some constant $C>0$. As $u \in L^{6}(\R^2)$, it follows that $h(u) \in L^{6/(p+1)}(\R^2)$ for $p\neq  3.$  While for $p=3$
$$
\int_{\R^2} |u|^{q (p+1)}dx=\int_{\R^2} |u|^{4q}dx <\infty, \ \mbox{for all}\quad q \geq 1,
$$
in particular,  $V_{xx} h(u) \in L^{q}(\R^2)$ for $q \in (1,3/2)$. 

For the second term, $(h_2)$ leads to
$$
|h'(u) u_x|\leq C |u|^{p}|u_x|.
$$
Note that  $|u|^{p}|u_x|$ belongs to $ L^{6/(p+1)}(\R^2).$ In fact, from (\ref{Passo1b}) since $|u|_{L^{\infty}} \leq C,$ for $p \neq 3,$ then
$$
\int_{\R^2} |u|^{p 6/(p+1)}|u_x|^{6/(p+1)}dx\leq C(\int_{\R^2} |u_x|^{12/(p+1)}dx)^{1/2} <\infty, \ \mbox{if}\  p \neq 3.
$$
For $p=3,$ $ u \in L^{q}(\R^2) $ for all $ q> 1.$ Then 
$$
\int_{\R^2} |u|^{3q} |u_x|^q  dx\leq (\int_{\R^2} |u|^{3sq}dx)^{s}(\int_{\R^2} |u_x|^{s'q}dx)^{s'} <\infty
$$
where $s'q=12/(p+1)=3$. Then, $s'=3/q>1 \Leftrightarrow q<3.$

For third term 
$$
|h''(u)( u_x)^2|\leq C |u|^{p-1}|u_x|^2
$$ 
for some constant. We claim that  $|u|^{p-1}|u_x|^2$ belongs to $ L^{6/(p+1)}(\R^2).$ Indeed, from (\ref{Passo1b}) since $u \in L^{\infty}(\R^2)$ for $p \neq 3,$ then
$$
\int_{\R^2} |u|^{(p-1) 6/(p+1)}|u_x|^{12/(p+1)}dx\leq C(\int_{\R^2} |u_x|^{12/(p+1)}dx) <\infty.
$$
For $p=3,$ we know that $ u \in L^{q}(\R^2), $ for all $ q> 1.$ Then 
$$
\int_{\R^2} |u|^{q 2} |u_x|^{2q}  dx\leq (\int_{\R^2} |u|^{sq2}dx)^{s}(\int_{\R^2} |u_x|^{s'q2}dx)^{s'} <\infty,
$$
where $s'q2=12/(p+1)=3$, that is, $s'=3/2q>1$. Thereby, $s'>1 \Leftrightarrow q<3/2.$

Finally, the last term verifies 
$$
|h'(u) u_{xx}|\leq C |u|^p |u_{xx}|
$$ 
for some $C>0$. we claim that $ |u|^p |u_{xx}| \in L^{6/(p+1)}(\R^2).$ In fact, from (\ref{Passo1b}), $u \in L^{\infty}(\R^2)$ for $p \neq 3,$ then
$$
\int_{\R^2} |u|^{p 6/(p+1)}|u_{xx}|^{6/(p+1)}dx\leq C(\int_{\R^2} |u_{xx}|^{6/(p+1)}dx) <\infty, \ \mbox{if}\ p \neq 3.
$$
For $p=3,$ $ u \in L^{q}(\R^2)$ for all $ q> 1.$ Thus,
$$
\int_{\R^2} |u|^{q 3} |u_{xx}|^{q}  dx\leq (\int_{\R^2} |u|^{sq3}dx)^{s}(\int_{\R^2} |u_{xx}|^{s'q}dx)^{s'} <\infty,
$$
where $s'q=6/(p+1)=6/4$, that is, $s'=3/2q>1$. Note that, $s'>1 \Leftrightarrow q<3/2.$ \\

The above analysis proves the Claim.

\vspace{0.5cm}

Differentiating (\ref{eq9}) with respect to $x$, twice, we get

\begin{equation}\label{eq22}
  -\Delta ( u_{xx}) + (u_{xx})_{xxxx} =  f_{xx} , \quad  (x,y)\in \R \times \R.
\end{equation}
Applying the Fourier transfom in the equation (\ref{eq22}), as in (\ref{eq10}), we have
$$
\widehat{u_{xx}}(\xi)=-(\frac{|\xi_1|^2}{|\xi|^2+|\xi_1|^4})\widehat{f}\equiv -\Phi_1(\xi)\widehat{f},
$$
or equivalently,
\begin{equation}\label{eq23} 
u_{xx}=\overbrace{(-\Phi_1(\xi)\widehat{f})}^{\vee}.
\end{equation}

By \cite[Corollary 1]{lizorkin}, $\Phi_1$ is a Fourier multiplier on $L^q(\R^2)$ for $ q =6/(p+1)$ if  $p \neq  3,$ and $1\leq q < 3/2$ if $p=3.$ Therefore,
\begin{equation}\label{regularuxx1}
\begin{array}{rl}
 u_{xx} \in L^{6/(p+1)}(\R^2),& \quad \mbox{if} \quad p \neq  3\\
u_{xx} \in L^{q}(\R^2), \ \forall q, \ 1\leq q< 3/2,& \ \quad \mbox{if} \quad  p=3
\end{array}
\end{equation}
Similarly, 
\begin{equation}\label{regularuxx11}
\begin{array}{rl}
 u_{xxxx}, u_{xxy},u_{yy} \in L^{6/(p+1)}(\R^2),&\quad \mbox{if} \quad  p \neq  3\\
 u_{xxxx}, u_{xxy},u_{yy}  \in L^{q}(\R^2), \ \forall q \in [1, 3/2), & \quad \mbox{if} \quad p=3. \\
\end{array}
\end{equation}

\vspace{0.5cm}

\noindent By (\ref{eq9}), 
\begin{equation}\label{eq121}
  -\Delta u + u =  \tilde{g}  \quad  \mbox{in} \quad \R^2 
\end{equation}
where
$$
\tilde{g}=f+u-u_{xxxx}.
$$
Since 
\begin{equation}\label{Passo1FFF}\left\{
\begin{array}{rl}
\tilde{g} \in L^{6/(p+1)}(\R^2),& \quad \mbox{if} \quad p\neq  3\\
\tilde{g} \in L^{q}(\R^2), \ \forall q \in [1, 3/2), & \quad \mbox{if} \quad p=3.
\end{array} \right.
\end{equation}
By \cite[Theorem 1]{Kavian},  
\begin{equation}\label{Passo1FFG}
u \in W^{2,q}(\R^2) \ for \ \left\{\begin{array}{rcl}
q = 6/(p+1)&\mbox{if}& p\neq  3\\
q\in (1, 3/2) & \mbox{if}& p=3,
\end{array} \right.
\end{equation}
Recalling that for any $q>1$the embedding
$$
W^{2,q}(\R^2)\hookrightarrow C^{0,\alpha}(\overline{\Omega})
$$
is continuous for any smooth bounded domain $\Omega \subset \R^2$ and $0 < \alpha \leq 2 -2/r,$ it follows that $ u \in C(\R^2).$ Moreover, by using bootstrapping arguments, there are $0<r_1<r_2$ and $C>0$ such that
$$
\|u\|_{W^{2,q}(B_{r_1}(x))} \leq C\|\tilde{g}\|_{L^{q}(B_{r_2}(x))}, \quad \forall x \in \R^2.
$$
Using Sobolev embeddings, there is $K>0$ independent of $x$ such that
$$
\|u\|_{C(\overline{B}_{r_1}(x))} \leq K\|\tilde{g}\|_{L^{q}(B_{r_2}(x))}, \quad \forall x \in \R^2.
$$
The last inequality gives
\begin{equation} \label{L1}
u(x) \to 0 \quad \mbox{as} \quad |x| \to +\infty. 
\end{equation}
This completes the proof of the theorem.  

\begin{flushright}$\Box$
\end{flushright}

\subsection{Concentration of the solution}

In this section, we study the concentration of the maximum point of the solution obtained in the previous section closed to the set where $V$ assumes its global maximum. In what follows, $u_\epsilon \in X$ denotes  a ground state solution obtained in Theorem \ref{teorema1}, that is,
$$
I_\epsilon (u_\epsilon)=c_\epsilon \quad \mbox{and} \quad I'_\epsilon (u_\epsilon)=0.
$$
In addition, we will fix $\epsilon_n \rightarrow 0,$ $u_n:=u_{\epsilon_n},$$c_n:=c_{\epsilon_n}$ and $I_n:=I_{\epsilon_n}.$

\begin{lem}\label{lema3} There exist constants $R,\eta>0$ and a sequence $(y_n)\subset \R^2$ such that
$$\lim_{n \rightarrow \infty} \int_{B_{R}(y_n)}|u_n|^2 dx dy \geq \eta.$$
\end{lem}

\noindent {\bf Proof.} If the lemma does not occur, by applying the Lions' Lemma version for $X$ found in \cite{Willem}, we have the limit 
$$
u_n \rightarrow 0 \quad \mbox{in}\quad L^4(\R^2),
$$
which leads to 
$$
u_n \rightarrow 0 \quad \mbox{in}\quad X.
$$
Therefore
$$
	c_n=I_n (u_n)=\frac{1}{2}\|u_n\|^2- \frac{1}{4}\int_{\R^2}V(\epsilon_n (x,y) H(u_n)  dx dy=o_n(1)
$$
which is a contradiction, because 
$$
\lim_{n \rightarrow \infty}c_n=c_0>0.
$$

\begin{flushright}$\Box$
\end{flushright}

From now on, set 
$$
w_n(x,y)=u_n((x,y)+y_n).
$$
After changing variable, we see that
\begin{equation}\label{15}
(w_n,v)=\int_{\R^2}V(\epsilon_n(x,y)+\epsilon_n  y_n)h(w_n) v  dx dy,\ \ \forall v \in X \quad \mbox{and} \quad n \in \mathbb{N}
\end{equation}
Once $ \|w_n\|=\|u_n\|$ and $(u_n)$  is bounded in $X,$  the sequence $(w_n)$ is bounded in $X.$ Consequently, up to a subsequence,
$$
w_n \rightharpoonup w \quad \mbox{in }\  X,
$$

$$
w_n \rightarrow w \quad \mbox{in }\  L^{2}_{loc}(\R^2)
$$
and
$$
w_n \rightarrow w \quad \mbox{a.e. in }\  \R^2.
$$
From the above limits, 
$$
\int_{ B_{R}(0)}|w_n|^2 dx dy=\int_{ B_{R}(y_n)}|w_n|^2 dx dy\geq \eta>0,
$$
and so, 
$$
\int_{ B_{R}(0)}|w|^2 dx dy\geq \eta,
$$
from where it follows $w \neq 0.$

\begin{lem}\label{lema4} For some subsquence of $\{w_n\}$, still denoted by itself, we have 
$$
w_n \rightarrow w \quad \mbox{in }\  X
$$ 
and
$$
\epsilon_n y_n\rightarrow y^* \in \mathcal{V}=\{z \in \R^2: V(z)=V_0\}.
$$
\end{lem}

\noindent{\bf Proof.} Arguing as in the proof of Lemma \ref{lema1},  we infer that
$$
\|w\|^2\leq \int_{\R^2}\limsup_{n \to +\infty}  V(\epsilon_n(x,y)+ \epsilon_n y_n)H(w)  dx dy.
$$

The above inequality permits to prove the following claim
\begin{claim} \label{CL2}
$(\epsilon_n y_n)$ is bounded in $\R^2.$
\end{claim}

\noindent{\bf Verification.} Suppose by contradiction that $(\epsilon_n y_n)$  possesses a subsequence, still denoted by $(\epsilon_n y_n)$, such that $|\epsilon_n y_n|\longrightarrow +\infty.$ By (\ref{Vinfty}), 
\begin{equation}\label{66} \int_{\R^2}\limsup_{n\to +\infty}  V(\epsilon_n (x,y)+ \epsilon_n y_n) H(w)dx dy\leq  \int_{\R^2}  V_{\infty} H(w) dx dy.
\end{equation}
Therefore, from (\ref{66}), 
$$
\|w\|^2\leq  \int_{\R^2}  V_{\infty} H(w) dx dy,
$$
or equivalently 
\begin{equation}\label{77} 
I'_{\infty}(w)w\leq 0.
\end{equation}
Let $t>0 $ be a number such that
$$
I_{\infty}(tw)=\max_{s \geq 0}I_{\infty}(sw).
$$
From (\ref{77}), we can guarantee that $t \in (0,1].$ Arguing as (\ref{8}),  
$$
c_\infty \leq \liminf_{n \to +\infty} I_n(u_n)=\liminf_{n \to +\infty}c_n=c_0,
$$
which is absurd. Therefore, we can suppose that, up to a subsequence, 
$$
\epsilon_n y_n \to y^*,
$$
for some $y^*\in \R^2$. Now, using the last limit, we are able to show the following claim
\begin{claim} \label{CL3}  
$$
y^* \in \mathcal{V}=\{z \in \R^2: V(z)=V_0\}.
$$
\end{claim} 
\noindent{\bf Verification.} If $y^* \notin \mathcal{V}$, we must have
$$
V(y^*)< V_0,
$$
and so, 
\begin{equation}\label{16}c_{V(y^*)}>c_{V_0}=c_0.
\end{equation}
Now repeating the arguments made above, just changing $c_\infty$ 
by $c_{V(y)},$ and  $ I_\infty$ by $I_{V(y^*)},$ we obtain
\begin{eqnarray*}
c_{V(y^*)}&\leq& I_{V(y^*)}(tw)=I_{V(y^*)}(tw)- \frac{1}{2}I'_{{V(y^*)}}(tw)tw\\
&\leq&\liminf_{n \to +\infty}c_n=c_0,
\end{eqnarray*}
which contradicts (\ref{16}).

Now, with the behavior of $(\epsilon_n y_n)$ in our hands, we can prove that $(w_n)$ converges strongly to $w$ in $X$.

\begin{claim} \label{CL4} 
$$
w_n \to  w \quad \mbox{in }\  X.
$$
\end{claim}

\noindent{\bf Verification.} Notice that by (\ref{15}), 
$$
I'_0(w)=0 \quad \mbox{and} \quad  w\neq 0.
$$
Therefore,
\begin{eqnarray*}
c_0&\leq& I_0(w)=I_0(w)- \frac{1}{2}I'_{0}(w)w=\int_{\R^{2}}\frac{V(y^*)}{2}(h(w)w-2H(w))\\
&\leq& \liminf_{n \to +\infty}\int_{\R^{2}}\frac{V(\epsilon_n x, \epsilon_n y+\epsilon_n y_n)}{2}(h(w_n)w_n-2H(w_n))\\
&\leq&\limsup_{n \to +\infty}\int_{\R^{2}}\frac{V(\epsilon_n x, \epsilon_n y+\epsilon_n y_n)}{2}(h(w_n)w_n-2H(w_n))\\
&\leq&\limsup_{n \to +\infty}\int_{\R^{2}}\frac{V(\epsilon_n x, \epsilon_n y)}{2}(h(u_n)u_n-2H(u_n))\\
&=&\lim_{n \to +\infty}( I_n(u)- \frac{1}{2}I'_{n}(u_n)u_n)\\
&=&\lim_{n \to +\infty} I_n(u_n)=\lim_{n \to +\infty}c_n=c_0.
\end{eqnarray*}
From this, 
$$
\lim_{n \to +\infty}\int_{\R^{2}}\frac{V(\epsilon_n x, \epsilon_n y+\epsilon_n y_n)}{2}(h(w_n)w_n-2H(w_n))=\int_{\R^{2}}\frac{V(y^*)}{2}(h(w)w-2H(w).
$$
The above limit combined with $(h_3)$ give
$$
V(\epsilon_n x, \epsilon_n y+\epsilon_n y_n)h(w_n)w_n \to V(y^*)h(w)w \quad \mbox{in} \quad L^{1}(\R^2)
$$
and
$$
V(\epsilon_n x, \epsilon_n y+\epsilon_n y_n)H(w_n) \to V(y^*)H(w) \quad \mbox{in} \quad L^{1}(\R^2).
$$
Since 
$$
\|w_n\|^{2}=\int_{\R^2}V(\epsilon_n x, \epsilon_n y+\epsilon_n y_n)h(w_n)w_n
$$
and
$$
\|w\|^{2}=\int_{\R^2}V(y^*)h(w)w,
$$
we deduce that
$$
\lim_{n\rightarrow \infty}\|w_n\|^2=\|w\|^2.
$$
As $w_n \rightharpoonup w \ \mbox{in} \ X$ and $X$ is a Hilbert space, we conclude
$$
w_n \to w \ \mbox{in} \ X.
$$
This proves the Lemma.

\begin{flushright}$\Box$
\end{flushright}

In the sequel, we set
$$
g_n(x,y)=-\frac{\partial^{2}}{\partial x^{2}}(V(\epsilon_n(x,y)+\epsilon_n y_n)) h(w_n)- \frac{\partial^{4}}{\partial x^{4}}w_n+w_n.
$$
Hence, $w_n$ satisfies the following equation
$$
-\Delta{w_n}+w_n=g_n \quad \mbox{in} \quad \mathbb{R}^2.
$$
Repeating the same arguments explored in the previous section, we have $w_n \in W^{q}(\R^2)$ where 
\begin{equation}\label{Passo1FFG'}
\left\{
\begin{array}{rl}
q=6/(p+1),&  \quad \mbox{if} \quad p\neq  3\\
q, & \quad \mbox{for} \quad q \in [1,3/2), \quad \mbox{if} \quad  p=3
\end{array} \right.
\end{equation}
Once $w_n \to w$ in $X$, the multiplier Fourier used in Section 3 permits to prove that 
\begin{equation} \label{G1}
g_n \to \hat{g} \quad \mbox{in}  \quad L^{q}(\R^2),
\end{equation}
where	
$$
\hat{g}(x,t)=-V(y^{*})\frac{\partial^{2}}{\partial x^{2}}(h(w))(x,t)- \frac{\partial^{4}}{\partial x^{4}}w(x,t)+w(x,t).
$$
Then, by using again bootstrapping arguments,  there are $0<r_1<r_2$ and $C>0$ such that
$$
\|w_n\|_{W^{2,q}(B_{r_1}(x))} \leq C\|g_n\|_{L^{q}(B_{r_2}(x))}, \quad \forall x \in \R^2 \quad \mbox{and} \quad n \in \mathbb{N}.
$$
Using Sobolev embeddings, there is $K>0$ independent of $n$ and $x$ such that
$$
\|w_n\|_{C(\overline{B}_{r_1}(x))} \leq K\|g_n\|_{L^{q}(B_{r_2}(x))}, \quad \forall x \in \R^2 \quad \mbox{and} \quad n \in \mathbb{N}.
$$
The last inequality together with (\ref{G1}) gives
\begin{equation} \label{L1}
w_n(x) \to 0 \quad \mbox{as} \quad |x| \to +\infty \quad \mbox{uniformly in} \quad n.
\end{equation}
Since $\|w_n\|_{\infty}=\|u_n\|_{\infty}$ and $\|u_n\|_\infty \not\to 0$, we derive that there is $\delta_0>0$ such that
\begin{equation} \label{L2}
\|w_n\|_{\infty} \geq \delta_0, \quad \forall n \in \mathbb{N}. 
\end{equation}
From (\ref{L1}) and (\ref{L2}), there exists $R>0$ such that
$$
\max_{z \in\R^2}|w_{n}(z)|=\max_{z \in \overline{B}_R(0)}|w_{n}(z)|.
$$ 
In what follows, we denote by $z_n \in \overline{B}_R(0)$ a global maximum point of $|w_n|$. Thereby, $\xi_n=z_n+y_n$ is a global maximum point of $u_n$, that is,
$$
|u(\xi_n)|=\max_{z \in\R^2}|u_{n}(z)|.
$$
Recalling that $(z_n)$ is bounded and $\epsilon_n y_n \to y^{*} \in \mathcal{V}$, we deduce that
$$
\lim_{n \to +\infty}V(\epsilon_n \xi_n)=V(y^{*})=V_0.
$$
This proves the concentration phenomenon.

\section{Final remarks}
After concluding this paper, the authors have observed that is possible to prove a multiplicity result of solutions by using Lusternik-Schnirelmann category. To this end, it is enough to adapt the arguments found in Alves and Figueiredo \cite{AG2006}. The result of multiplicity can be state of the following way: 
\begin{thm} For each $\delta >0$ small enough, there is $\epsilon^*>0$ such that problem (\ref{eq6}) has at last $cat_{\mathcal{V}_\delta}(\mathcal{V})$, where 
	$$
	\mathcal{V}_\delta=\{z \in \R^2 \,:\, dist(z,\mathcal{V})\leq \delta \}
	$$  
	for all $\epsilon \in (0,\epsilon^*)$. Moreover, if $u_\epsilon$ denotes one of these solutions and $q_\epsilon$ is a global maximum point of $|u_\epsilon|$, we have that
	$$
	\lim_{\epsilon \to 0}V(\epsilon q_\epsilon)=\max_{z \in \R^2}V(z).
	$$	
\end{thm}  
We would like to point out that if $Y$ is a closed subset of a topological space $X$, the Lusternik-Schnirelman category $cat_{X}(Y)$ is the least number of closed and contractible sets in $X$ which cover $Y$. If $X=Y$, we use the notation $cat(X)$. For more details about the Lusternik-Schnirelman category see \cite{Willem}.




\end{document}